\documentclass[12pt]{{amsart}}
\usepackage[a4paper,top=3cm,bottom=2cm,left=2cm,right=2cm,marginparwidth=1.75cm]{geometry}

\usepackage[utf8]{inputenc}
\usepackage{amsmath}
\usepackage{amsfonts}
\usepackage{amssymb}
\usepackage{graphicx}
\usepackage[linesnumbered,ruled]{algorithm2e}
\usepackage[T1]{fontenc}
\usepackage{float}
\usepackage{tikz}
\usepackage{color}
\usepackage{verbatim}
\usepackage{listings}
\usepackage{graphics,graphicx}                    
\usepackage{pstricks,pst-node,pst-tree}
\usepackage{pstricks,pst-node,pst-tree,pstricks-add} 
\usepackage[all]{xy}
\usepackage{amsmath}
\usepackage{graphicx}
\usepackage{url}
\usepackage{amsfonts}
\usepackage{amssymb}
\usepackage{hyperref}
\usepackage{lastpage} 
\usepackage{fancyhdr}
\usepackage[colorinlistoftodos]{todonotes}
\newtheorem{Theorem}{Theorem}[section]

\numberwithin{equation}{subsection}
\newtheorem{Proposition}{Proposition}[section]
\newtheorem{Lemma}{Lemma}[section]
\newtheorem{Corollary}{Corollary}[section]
\newtheorem{Remark}{Remark}[section]

\title{\small\textbf{Rational and isolated qaudratic points on hyperelliptic curves of rank 0 and small genus}}
\author [Brice M.Miayoka]{Brice M.  Miayoka}
\address{ Brice Moussolo Miayoka\\ 
Department of Mathematics \\ 
Marien NGOUABI University   \\ 
Brazzaville Congo}
\email{bricemiayo@gmail.com}

\begin{document}
\maketitle 
\begin{abstract}
In this article, we present a method for computing rational points on hyperelliptic curves of genus~3 and isolated quadratic points on hyperelliptic curves of genus~2 and~3 whose Jacobians have rank~0. Our approach begins by computing the image of the Mordell--Weil group on the associated Kummer variety and then determining which of these points correspond to rational or isolated quadratic points on the curve. We have developed and implemented this algorithm using the computer algebra system \texttt{\texttt{Magma}}. The method takes advantage of structural properties specific to hyperelliptic curves and their Jacobians. We applied our algorithm to a large dataset, analyzing 7,396 genus~3 hyperelliptic curves and 12,075 genus~2 hyperelliptic curves.

\end{abstract}
\begin{center}
\section{Introduction}
\end{center}
Let \( C \) be a hyperelliptic curve of genus \( g\geq 2 \) defined over \( \mathbb{Q} \).  While it is known by Faltings's theorem \cite{FG} that the set \( C(\mathbb{Q}) \) is finite, computing it is generally  a difficult problem. However, when the rank \( r \) of the Jacobian \( J_C \) of \( C \) over $\mathbb{Q}$ is low, certain techniques make it possible. For example, when $r< g$, one can use the method of Chabauty and Coleman \cite{cl, CR}.  When $g=3$, then we have the algorithm of
Balakrishnan, Bianchi, Cantoral-Farfán, Çiperiani, and Etropolski, which computes \( C(\mathbb{Q}) \) when \( r = 1 \) \cite{BBC} and the algorithm of Fernández and Hashimoto \cite{Fu12} which computes \( C(\mathbb{Q}) \) for rank \( r = 0 \). However, these algorithms rely heavily on  algorithms for  Coleman integrals \cite{Bjen12} and they are currently restricted to curves with a rational Weierstrass point.

On the other hand, Cassels and Flynn, in  \cite[Section 13.1]{Ca12}, developed a very efficient method for computing \( C(\mathbb{Q}) \) in the case where \( g = 2 \) and \( r = 0 \). This method is simple and is not based on the computation of Coleman integrals. This algorithm was implemented in \texttt{Magma} by Stoll \cite{st} and can be called in \texttt{Magma} as \texttt{Chabauty0}.
In genus 3, such a technique has not yet been described, which is why in this paper, we introduce an algorithm to compute rational points on hyperelliptic curves of genus 3 with a rank 0 Jacobian which does not require the existence of a rational point on the curve. We have implemented and executed our algorithm on 7396 hyperelliptic curves of even and odd degree. We found 503 hyperelliptic curves with even models that have no rational points.

We also show that our method can do more than simply compute \( C(\mathbb{Q}) \): it also allows us to determine the isolated quadratic points  on hyperelliptic curves of genus 2 and 3 whose Jacobian \( J_C \) has rank 0.
Let \( k \) be a finite extension of \( \mathbb{Q} \) of degree \( d > 0 \), and let \( P = (x, y) \in C(k) \) be a point defined over \( k \). Three cases can arise concerning the nature of the point \( P \):
\begin{itemize}
  \item \( x \in \mathbb{Q} \): in this case, the point \( P \) is invariant under the action of the Galois group \( G := \operatorname{Gal}(k/\mathbb{Q}) \) on the \( x \)-coordinate, and it arises from the natural covering \( C \to \mathbb{P}^1 \). Points of this type are generally infinite in number.

  \item The point \( P \in C(k) \) comes from a covering \( C \to E \), where \( E \) is an elliptic curve defined over \( \mathbb{Q} \) such that \( E(\mathbb{Q}) \) has strictly positive rank. This yields infinitely many points of degree \( d \). This case is excluded here, because such a factorization would imply that the Jacobian \( J_C \) does not have rank 0, contradicting our assumption.

  \item The remaining points \( P \in C(k) \) are said to be isolated, as they are neither defined over \( \mathbb{Q} \), nor do they arise from a covering of a curve with infinitely many rational points.
\end{itemize}
A theorem by Bourdon, Ejder, Liu, Odumodu, and Viray states that, under the above assumptions, the set $\mathcal{L}$ of the  \( k\)-rational points on \( C \) that do not arise from the natural covering \( C \to \mathbb{P}^1 \) are finite \cite{BEV}: these are precisely the isolated points. In general, computing such points is not easy. When the rank of the Mordell--Weil group \(J_C(\mathbb{Q})\) satisfies \( r \leq g - d \), Siksek~\cite{Sisk}, building on ideas of Klassen~\cite{Kla9}, introduced a method known as symmetric Chabauty to compute these points. Our approach here is different from his approach.  Our algorithms relies on the algorithm developed by Müller and Reitsma, described in \cite{Mur13}, which computes \( J_C(\mathbb{Q})_{\text{tors}} \).


This paper is structured as follows: in Section~\ref{sec:2}, we describe the computation of rational points on curves of genus $g \geq 2$ whose Jacobian has rank 0. In Section~\ref{sec:3}, we recall the algorithm described by Müller and Reitsma for computing the structure of $J_C(\mathbb{Q})_{\text{tors}}$, as well as its image on the Kummer variety of the Jacobian of a  hyperelliptic curves of genus 3. We describe our algorithm  to compute rational points on hyperelliptic curves in Section~\ref{sec:4}. In Section \ref{sec:51}, we present the computation of degree 2 isolated points on hyperelliptic curves of genus 2 and 3. In Section \ref{sec:5}, we compare it with the approach of Fern\'andez and Hashimoto,  and mention some possible generalizations. In Section~\ref{sec:6}, we analyze our data and present some examples.
\section*{Acknowledgments.} This study was carried out with the financial support of the Coimbra Group, through the Scholarship Programme for Young African Researchers. It is with great pleasure that I thank Steffen Müller for suggesting this research topic, for answering our numerous questions, as well as for his valuable suggestions and insightful comments on the preliminary versions of this article.
I also express our deep gratitude to the Boston University Library (USA) and to the Artificial Intelligence Group at the Bernoulli Institute of the University of Groningen (Netherlands), who generously provided us with access to their servers, which were used for our computations.
I am especially grateful to Jennifer Balakrishnan for sponsoring my affiliation with Boston University (USA).
Finally, I would like to express my sincere thanks to my supervisors, Tony Ezome and Régis Babindamana, for their continuous support and encouragement throughout this work.

\begin{center}
\section{Computing  the rational points on a curve when its Jacobian has  rank 0}\label{sec:2}
\end{center}
Let \( C \) be a  smooth projective geometrically irreducible curve of genus \( g \geq 2 \), defined over \( \mathbb{Q} \). We denote by \( J_C \) the Jacobian of  \( C \). We assume  throughout that \( J_C \) has rank 0 over \(  \mathbb{Q} \). In this section, we discuss a general strategy to compute the set  \( C(\mathbb{Q}) \).
\\ Assume we know some   \(P_0\in C(\mathbb{Q}) \), so there exists an Abel-Jacobi embedding \( \iota \colon C \to J_C \) defined over \( \mathbb{Q} \), which maps a point \( Q \in C \) to  the divisor class \( [(Q) - (P_0)] \). Then the image \( \iota(C(\mathbb{Q})) \) will consist of $\mathbb{Q}$-rational  divisor classes. This image is generally a finite subset of \( J_C(\mathbb{Q}) = J_C(\mathbb{Q})_{\text{tors}} \), and \( C(\mathbb{Q}) \) can be computed as follows.
  \begin{enumerate}
    \item[(a)] \textit{Computing the torsion subgroup of \( J_C(\mathbb{Q}) \)}. This is actually the whole group of rational points, since the rank is 0.
    \item[(b)] \textit{Recognizing the image of \( C(\mathbb{Q}) \) inside \( J_C(\mathbb{Q}) \) under \( \iota \)}.
\end{enumerate}
If \( P_0 \) is a Weierstrass point, we can experiment with the approach developed by Poonen in \cite{zbMATH01714029}.\\
In the case where \( C \) is a hyperelliptic curve of genus \( g \geq 2 \), , it can be defined by an equation  \( C: y^2 = f(x) \), where \( f \) is a polynomial with coefficients in \( \mathbb{Q} \) without repeated roots and \( \deg(f) = 2g + 2 \) or \( 2g + 1 \). Then we can instead use a slightly different approach, as we now explain. It has the advantage that we do not need a rational point $P_0\in C(\mathbb{Q})$ to begin with.\\
 Let \( \phi \colon C \to C \) be the hyperelliptic involution. The divisor \( P + \phi(P) \) is always \(\mathbb{Q}\)-rational for any \( P \in C \), so the class \( c \) of this divisor is independent of \( P \).
We therefore define a map \( \lambda \colon C \to J_C \) that sends \( P \) to \( [2(P)] - c \). We use
\( \{(x_1, y_1), (x_2, y_2)\} \) as shorthand notation for the divisor class containing \( [(P_1) + (P_2)]-c \), where \( P_1 = (x_1, y_1) \) and \( P_2 = (x_2, y_2) \in C( \mathbb{Q})\). We seek to find the points in \( J_C(\mathbb{Q}) \) that are of the form \( \{P, P\} \).
\begin{Lemma}\label{zz}
 Suppose $(x_2, y_2) \ne (x_1, -y_1)$. Then the divisor class \( \{(x_1, y_1), (x_2, y_2)\} \) is of the form \( \{P, P\} \) if and only if \( (x_1 + x_2)^2 - 4 x_1 x_2 = 0 \).
\end{Lemma}
Applying this result in practices requires knowledge of \( J_C(\mathbb{Q})_{\text{tors}} \). An algorithm to compute the structure of \( J_C(\mathbb{Q})_{\text{tors}} \) for genus 3 hyperelliptic curves over \( \mathbb{Q} \) was proposed by Müller and Reitsma in \cite{Mur13}.
We recall their algorithm in the following section.
\begin{center}
 \section{Computing torsion subgroups of Jacobians of genus 3 hyperelliptic curves}\label{tors}\label{sec:3}
\end{center}

\subsection{Kummer variety\label{8}}
Let \( C \) be a hyperelliptic curve of genus \( g \) defined over \( \mathbb{Q} \), and let \( J_C \) be its Jacobian. The associated Kummer variety, denoted \( K \), is the quotient of \( J_C \) by the involution \( P \mapsto -P \). This variety can be embedded into a projective space of dimension \( 2^g - 1 \) \cite{Mur12}.
\noindent In particular, when \( g = 2 \), the Jacobian is an abelian surface, and the associated Kummer variety is a quartic surface in \( \mathbb{P}^3 \), possessing 16 rational singular points, corresponding to the 2-torsion points of the Jacobian \( J_C \). The class of a degree 2 divisor, \( D = (x_1, y_1) + (x_2, y_2) \), is represented in \( K \) by the projective coordinates \( (\sigma_1 : \sigma_2 : \sigma_3 : \beta_0) \), where $\sigma_1=1$, $\sigma_2=x_1+x_2$,  $\sigma_3=x_1x_2$ and \( \beta_0 \) is an invariant associated with the curve \( C \) \cite[Section 3.1]{Ca12}. Cassels and Flynn in \text{ \cite[Section 13.1]{Ca12}} use this representation to compute rational points on \( C \). They show that these points satisfy the equation:
$(x_1 + x_2)^2 - 4x_1 x_2 = 0$

\noindent In genus 3, let us first recall that a divisor \( D \in \operatorname{Div}(C) \) is said to be in  \textit{general position} if it is effective and there does not exist a point in \( \operatorname{Supp}(D) \) such that
$D \geq (P) + (\phi(P))$ where \( P \in C \). In \cite[\S 2]{st12}, Stoll shows  the following proposition:
\begin{Proposition}\label{3}  Suppose $C: y^2=f(x)$. Let \( D_{\infty} \in \operatorname{Div}(C) \) be a divisor defined by:
$D_{\infty} = 2(\infty) \quad \text{if }  \text{deg}(f)= 7,$ where $D_{\infty} = (\infty_1) + (\infty_2) \quad \text{if }  \text{deg}(f)= 8.$  Then, for any nonzero point \( Q \) in \( J_C \), exactly one of the following conditions holds:
\begin{enumerate}
    \item[(a)] \( Q \) is the class of a divisor of the form \( D_Q - 2D_{\infty} \), where \( D_Q \) is a degree \( 4 \) divisor in general position.
    \item[(b)] \( Q \) is the class of a divisor of the form \( D_Q - D_{\infty} \), where \( D_Q \) is a degree \( 2 \) divisor in general position.\\
    In case (b), the divisor \( D_Q \) is uniquely determined by  \( Q \).
\end{enumerate}
\end{Proposition}
In this paper, we consider only case (b) of Proposition \ref{3}, as the points in case (a) cannot contain any element in the image of \( C(\mathbb{Q}) \) under the map \( \lambda \), defined in Section \ref{sec:2}. Note that the divisor $D_{\infty}$ is a representative of the class c defined above. In \cite{st12}, Stoll used this representation to explicitly construct the Kummer variety of \( J_C \). Indeed, there exists a canonical theta divisor \( \Theta \) on \( J_C \), whose support consists of \( 0 \) and the classes \( Q \in J_C \) corresponding to the degree \( 2 \) divisors mentioned in the above proposition. Any basis \( (\kappa_1, \dots, \kappa_8) \) of the Riemann-Roch space \( L(J_C,2\Theta) \) defines a rational map $
\kappa: J_C \to \mathbb{P}^7
$ such that \( \kappa(J_C) \) gives the Kummer variety \( K \) of \( J_C \).
Suppose that the curve \( C/\mathbb{Q} \) is defined by the equation:
\begin{equation}\label{eq:1}
C: y^2 = f(x), \end{equation}
where $f \in \mathbb{Q}[x]$  is squarefree of degree $7$ or $8$.
Let \( Q \in J_C(\mathbb{Q}) \) be a point of the form: $
Q = [D_Q - D_{\infty}],$
where $D_Q = (x_1, y_1) + (x_2, y_2) $
with \( (x_i, y_i) \in C(\mathbb{Q}) \) for \( i = 1,2 \). We define the polynomial:

\[
A(x,z) = (z_1 x - x_1 z)(z_2 x - x_2 z) = a_0 x^2 + a_1 xz + a_2 z^2 \in \mathbb{Q}[x,z].
\]
Stoll describes an explicit map kappa in \cite[Theorem 2.5 ]{st12}. For the point $Q$, we have::
\[
\kappa(Q) = (0 : a_0^2 : a_0 a_1 : a_0 a_2 : a_1^2 - a_0 a_2 : a_1 a_2 : a_2^2 : \sigma_8).
\]

\noindent We can rewrite this in terms of the coordinates of $(x_1, y_1)$  and $(x_2, y_2)$ as follows.
$$Q=\{(x_1, y_1),(x_2, y_2)\}\longrightarrow \kappa(Q)=(\sigma_1:\sigma_2:\sigma_3:\sigma_4:\sigma_5:\sigma_6:\sigma_7:\sigma_8),
$$ where $$\sigma_1=0, \sigma_2=1,\sigma_3=-(x_1+x_2), \sigma_4=x_1x_2, \sigma_5=x^2_1+x_1x_2+x^2_2, \sigma_6=-(x_1+x_2)x_1x_2, \sigma_7=(x_1x_2)^2,$$  and  $$\sigma_8=\dfrac{2y_1y_2-G(x_1,x_2)}{(x_1-x_2)^2}$$ satisfying the equation \begin{equation}\label{ss}
((x_1-x_2)^2\sigma_8+G(x_1,x_2))^2-4f(x_1)f(x_2)=0
\end{equation}
where
$$G(x_1,x_2)=2\displaystyle \sum^{4}_{j=0}f_{2j}(x_1x_2)^{j}+(x_1+x_2) \displaystyle \sum^{3}_{j=0}f_{2j+1}(x_1x_2)^{j}.$$\\In the case where  $x_1= x_2$,
we have $\sigma_i$ with $i=1,..,7$ as before but  we can calculate $\sigma_8$ by writing  \ref{ss} as \begin{equation}s_2\sigma^2_8+s_1\sigma_8+s_0=0.\end{equation} Then $s_2=0$  and $s_1=-2G(x_1,x_1)=-4f(x_1)$. If $s_1\neq 0$, then $Q\neq 0$ and  it follows that  $$\sigma_8=\dfrac{-s_0}{s_1}$$
In addition to the description of $K$, we also need algorithms for multiplying points in $K$ by a scalar, as well as for the sum and difference of two points in $K$. We can perform these operations using the following two maps: \begin{itemize}
                                                                                                                                                                                                                                                                                                   \item The map $\text{[[2]]} : K \rightarrow K$ such that $\kappa(2Q) = \text{[[2]]}(\kappa(Q))$ for all $Q \in J_C$. Let $R = \kappa(Q)$. Then $\text{[[2]]}(R) = (\delta_1(R) : \ldots : \delta_8(R))$, where $\delta_1, \ldots, \delta_8 \in \mathbb{Z}[f_0, \ldots, f_8][\sigma_1, \ldots, \sigma_8]$ are homogeneous quartic polynomials satisfying $(\delta_1, \delta_2, \ldots, \delta_8)(0, 0, \ldots, 0, 1) = (0, 0, \ldots, 0, 1)$. For more details on the construction of these polynomials, see \cite[Theorem 7.3]{st12}.
\item The map $B : \text{Sym}^2(K) \rightarrow \text{Sym}^2(K)$ such that for all $Q_1, Q_2 \in J_C$, we have
\[
B(\{\kappa(Q_1), \kappa(Q_2)\}) = \{\kappa(Q_1 + Q_2), \kappa(Q_1 - Q_2)\}.
\]

 Formulas for $B$ are given in \cite[ Lemma 8.1]{st12}.                                                                                                                                                                                                                                                                                               \end{itemize}
 In  \cite[ Lemma 4.5]{Mur13}, Müller and Reistma provide conditions that allow one to determine whether a given point \( R \in K(\mathbb{Q}) \) lifts to \( Q \in J(\mathbb{Q}) \) under \( \kappa \), in the form of case (b) of Proposition   \ref{3}. We note that Stoll provides such conditions for case (a) \cite{st12}.
Since we are only interested in case (b), we provide, for our purposes, the following result, which is equivalent to  Lemma~\ref{zz} mentioned above.
\begin{Corollary}\label{z1}
Let \( R = (0 : \sigma_2 : \sigma_3 : \sigma_4 : \sigma_5 : \sigma_6 : \sigma_7 : \sigma_8) \in K(\mathbb{Q}) \). The preimage \( \kappa^{-1}(R) \) is of the form \( \{P, P\} \)  if and only if \( \sigma_2 \neq 0 \) and \( \sigma_3^2 - 4\sigma_2 \sigma_4 = 0 \).
\end{Corollary}
\subsection{Computing the torsion subgroup}\label{sec:31}
The Müller and Reitsma algorithm is an extension of the method developed by Stoll \cite{Stoll} for computing the rational torsion subgroup of Jacobians of  curves in genus $2$, specifically for the case of hyperelliptic curves of genus 3.\\
Let $C$ be a hyperelliptic curve of genus $ 3$ defined over $\mathbb{Q}$ by the equation \eqref{eq:1}. Let $J_C$ be the Jacobian of $C$. It is known from the Mordell-Weil theorem that
$$
J_C(\mathbb{Q}) \cong \mathbb{Z}^r \times J_C(\mathbb{Q})_{\text{tors}},
$$
where $r$ is the rank of $J$ and $J_C(\mathbb{Q})_{\text{tors}}$ is the finite torsion subgroup.

An upper bound for $J_C(\mathbb{Q})_{\text{tors}}$ can be computed by considering the injection from the rational torsion subgroup into the reduction of $J_C(\mathbb{Q})$ modulo an odd prime of good reduction. Indeed, let $p>2$ be a prime of good reduction for  $J_C$. We have the reduction homomorphism $\rho_p: J_C(\mathbb{Q}_p) \to J_C(\mathbb{F}_p)$. For an integer $m \geq 1$ such that $p \nmid m$, the restriction of the map $\rho$ to $J_C(\mathbb{Q}_p)[m]$ is injective  \cite[Theorem C.1.4]{zbMATH01466163}. If we take $S$ to be a set containing some primes of good reduction and compute $\# J_C(\mathbb{F}_p)$ for all $p \in S$, then $\# J_C(\mathbb{Q})_{\text{tors}} \mid \gcd_{p \in S} (\# J_C(\mathbb{F}_p))$.\\
The algorithm  is described as follows and has been implemented in \texttt{Magma}, see \cite{degond6}:
The algorithm starts by determining the bound \(\beta\), which represents the difference between the naive height and the canonical height. Indeed, let \( P = (x_1: \dots: x_8) \in \mathbb{P}^7(\mathbb{Q}) \). The naive exponential height of $P$  is defined as:
$H(P) = \max\left( |x_1|, \dots, |x_8| \right),$
where \( x_1, \dots, x_8 \) are integers that are coprime. The naive height on \( J_C(\mathbb{Q}) \) is the function $h : J_C(\mathbb{Q}) \to \mathbb{R}_{\geq 0}$, defined by
$
h(Q) := \log(H(\kappa(Q))),
$
where \( \kappa \) is the map introduced in the previous section. The function \( h \) is quadratic up to a bounded function, hence the canonical height is well-defined:
\[
\hat{h}(Q) = \lim_{n \to \infty} \frac{h(nQ)}{n^2},
\] The difference between the canonical height $\hat{h}$ and the naive height $h$ is known to be bounded on $J_C(\mathbb{Q})$. That is, there exists a real constant $\beta > 0$ such that
\[
|\hat{h}(Q) - h(Q)| < \beta \quad \text{for all } Q \in J_C(\mathbb{Q}).
\]
\noindent \textbf{Algorithm 1. Computing the Torsion Subgroup}

\noindent \textbf{Input:} A Jacobian \( J_C/\mathbb{Q} \) \\
\textbf{Output:} The invariant factors of \( J_C(\mathbb{Q})_{\text{tors}} \).

\begin{enumerate}
    \item Compute  an upper bound $\beta$ as in \cite[Section 4.7]{Mur13}.
    \item Compute a multiplicative upper bound \( t \) for the size of the torsion subgroup by determining the structure of \( J_C(\mathbb{F}_p) \) for a reasonable number of good odd primes \( p \).
    \item For each prime factor \( q \) of \( t \), compute the \( q \)-part of \( J_C(\mathbb{Q})_{\text{tors}} \) using \cite[Algorithm 3.14]{Mur13}
    \item Deduce the invariant factors of \( J_C(\mathbb{Q})_{\text{tors}} \) from the invariant factors of its \( q \)-parts.
\end{enumerate}
Stoll has developed in \cite[Corollary 10.3]{st12} a method to compute a suitable bound \( \beta \). The bound \( \beta \) ensures that the algorithm only considers rational points that are potential candidates for the torsion subgroup, thereby avoiding the exploration of an infinite or unnecessarily large search space. Step 2 of Algorithm 1 involves describing the structure of \( J_C(\mathbb{F}_p) \) for primes \( p \) of good reduction. We use the arithmetic on the Kummer variety \( K/\mathbb{F}_p \) and check whether points in \( K(\mathbb{F}_p) \) lift to \( J_C(\mathbb{F}_p) \).
  In practice, we work directly in \( J_C(\mathbb{F}_p) \) if the arithmetic in \( J_C(\mathbb{F}_p) \) is implemented, i.e, if \( C(\mathbb{F}_p) \) is nonempty. See \cite[Section 4.4]{Mur13}.
\subsection{Checking whether reduced points lift}\label{sec:32}
In this work, we will be led to compute the points \( R \) of \( K(\mathbb{Q}) \) that lift to a point \( Q \) of \( J_C(\mathbb{Q})_\text{tors} \). To do so, we proceed as follows. This technique was developed by Müller and Reitsma in \cite[Section 3]{Mur13} following the analogous algorithm of Stoll for  genus 2, see \cite[\S11]{Stoll}:

Let \( p \) be a prime of good reduction for \( J_C \). Let \( q \neq p \) be a  good prime \cite[Section 3]{Mur13} . For any \( P \in J_C(\mathbb{F}_p) \) of \( q \)-power order \( m \), we can determine the unique lift \( \tilde{\kappa}(P) \) in \( \kappa(J_C(\mathbb{Q}_p)[m]) \) to any desired precision \( p^N \).

Using \( \beta \), we fix \( N \) and construct a lattice \( L \) with the following property: if there exists a point \( R \in \kappa(J_C(\mathbb{Q}_p)[m]) \cap K(\mathbb{Q}) \) whose reduction modulo \( p^N \) corresponds to the approximation of \( \tilde{\kappa}(P) \), then the shortest non-trivial vector in \( L \) must be this point \( R \).
The existence of such a point can be checked using the LLL algorithm \cite{LLL}. If a candidate is found, it remains to check if it lifts to \( J_C(\mathbb{Q})[m] \). For odd values of \( q \), we refer to \cite[Algorithm 3.4]{Mur13}. In the special case where \( q = 2 \), the set \( J_C(\mathbb{Q})[2] \) can be computed efficiently by factoring \( f \) over \( \mathbb{Q} \) and over an explicit finite list of quadratic extensions, see \cite[Section 4.5]{Mur13}.

This algorithm works in a general setting; it can be applied even when the points on the Jacobian and their arithmetic are not implemented (i.e., the curve \( C \) does not have a known rational point). Note that when \( \deg(f) = 7 \), the algorithm of Müller and Reitsma allows us to obtain the generators of \( J_C(\mathbb{Q}) \). However, in general, it only provides the invariant factors. Nevertheless, from the code \cite{degond6}, one can extract the images of the generators on \( K \).

\section{Computing the rational points on a hyperelliptic curve of genus 3 when its Jacobian is irreducible of rank 0}\label{sec:4}
In this section we describe our agorithm to compute $C(\mathbb{Q})$. We have implemented this algorithm in \textbf{\texttt{Magma}}, and it is available at \href{https://github.com/Brice202145/Chabauty0Genus3}{https://github.com/Brice202145/Chabauty0Genus3}.\\ We require $C$ to be hyperelliptic of genus $3$ given by a equation \eqref{eq:1}. Morover, we assume that its Jacobian $J_C$  has Mordell-Weil rank $0$ over $\mathbb{Q}$, this can be checked in practice using 2-descent \cite{Rk}  .\\
We first note that it is straightforward to find all rational Weierstrass points since this only involves finding all roots of $f(x)$ in $\mathbb{Q}$. We denote by $S$ the set of all rational Weierstrass points of $C$. $$S:=\{(x,0)\in \mathbb{Q}^{2}\mid f(x)=0\}.$$
 Let $K = \kappa(J_C)$ be the Kummer variety of $J_C$. We will find the affine rational points of \( C \) that are not Weierstrass points. For that, we compute a certain \( TK \subseteq K(\mathbb{Q}) \).  We distinguish two cases:
\begin{enumerate}
 \item If $ \deg(f) =7$, we directly compute \( \kappa(J_C(\mathbb{Q})) \) for all \( Q \in J_C(\mathbb{Q})\). We then define:  $$TK:=\{R\in  K(\mathbb{Q}) \mid R=\kappa(Q) :\;Q\in J_C(\mathbb{Q}),\; Q\neq 0\}.$$
 \item If  \(\deg(f) =8 \),   the algorithm described by Müller and Reitsma (detailed in Subsection \ref{sec:31}),  computes the invariant factors of \( J_C(\mathbb{Q})_{\text{tors}} \) rather than the points it contains. This is because, when the curve has degree $8$, it is generally not possible to represent and add points on the Jacobian in \texttt{Magma}. Fortunately,  we do not need to know \( J_C(\mathbb{Q}) \) explicitly, but rather its image on the Kummer variety \( K \). Let \( p > 2 \) be a prime of good reduction for \( J_C \), and let $P\in J_C(\mathbb{F}_{p})$.  The procedure in Section \ref{sec:32} finds the unique point \( R \in \kappa(J_C(\mathbb{Q}_p)_{\text{tors}}) \cap \kappa(J_C(\mathbb{Q})) \) that reduces to \( \kappa(P) \), if such a point exists. Our ultimate goal is determine whether this  \(R\in K(\mathbb{Q}) $ lifts to a point in $J_C(\mathbb{Q})\) (without having to compute a lift of $R$ to $J_C$ explicitly). We directly use the function \textbf{doesLiftToJacobian} in the \texttt{Magma} code  of Müller and Reitsma \cite{degond6}. By iteration over all points $P\in J_C(\mathbb{F}_{p})$, we compute $\kappa(J_C(\mathbb{Q}))$.
 We set
 $$
TK := \left\{ R \in K(\mathbb{Q}) \;\middle|\;
\begin{array}{l}
R \in \kappa\left(J_C(\mathbb{Q}_p)_{\text{tors}}\right) \cap \kappa\left(J_C(\mathbb{Q})\right), \\
\text{et } R \text{ se réduit en } \kappa(P) \text{ pour tout } P \neq 0
\end{array}
\right\}.
$$

\end{enumerate}

\noindent According to Subsection \ref{8}, a point $R \in TK$ is of the form
$$R=(\sigma_1:\sigma_2:\sigma_3:\sigma_4:\sigma_5:\sigma_6:\sigma_7:\sigma_8)$$
We will determine the set of rational points \( P \in C(\mathbb{Q}) \) such that \( \{P, P\} = \kappa^{-1}(R) \). According to Corollary \ref{z1}, this occurs when

$$\sigma_1=0,\; \sigma_2\neq 0\; \text{and} \; \sigma^{2}_3-4\sigma_2\sigma_4=0$$
In this case the $x$-coordinate of $P$ is given by $x := \frac{-\sigma_{3}}{2 \sigma_{2}}$. and the $y$-coordinate is obtained by solving the equation $y^2=f(x)$. The indeterminacy of the sign of \( y \) does not affect the reasoning since we are only interested in rational points, and we always consider both points \( P = (x, y) \) and \( \phi(P) = (x, -y) \).  The points $P$  and $\phi(P)$ are added to $S$.\\
Finally, we will now add the points at infinity to \( S \). If our curve $C$ is of degree 7, this model has a unique rational point at infinity, which we denote \( \infty \). In the case where the curve $C$ is of degree 8, defined by the equation $C \colon y^2 := f_8 x^8 + f_7 x^7 + f_6 x^6 + f_5 x^5 + f_4 x^4 + f_3 x^3 + f_2 x^2 + f_1 x + f_0.$ It has two rational points at infinity, which we denote \( \infty_{1} \) and \( \infty_{2} \), provided that \( f_8 \) is a square. In this case, \( \infty_{1} = (1 : s : 0) \) and \( \infty_{2} = (1 : -s : 0) \), where \( f_8 = s^2 \) and \( s \in \mathbb{Q} \). Otherwise, there are not rational points at infinity.\\
The input of our algorithm consists of the hyperelliptic curve $C$.\\
The ouput is the set $S$ of rational points of $C$.
\begin{center}
\fbox{
\begin{minipage}{12cm}
\textbf{Algorithm 2: Rational points on the hyperelliptic curve \( C \colon y^2 = f(x) \), where \( f(x) \) is a squarefree polynomial of degree 7 or 8 and the Jacobian \( J_C(\mathbb{Q}) \) has rank 0.} \label{bb}\\
\textbf{Input:}  a genus 3 hyperelliptic curve  $C$.  \\
\textbf{Output:} The set $S$ of $\mathbb{Q}$-points.\\
1- Compute $S := \{(x,0) \mid x \in \mathbb{Q}, f(x) = 0\}$.\\
2- Compute the set $TK$ of points on the Kummer variety of $J$ defined over $\mathbb{Q}$.\\
3- \textbf{for } $R := (\sigma_1 : \sigma_2 : \sigma_3 : \sigma_4 : \sigma_5 : \sigma_6 : \sigma_7 : \sigma_8) \text{ in } TK$ \textbf{ do}\\ \hspace*{1cm} 4- \textbf{if } $\sigma_1=0$, $\sigma_{3}^2 = 4\sigma_{2}  \sigma_{4}$ \textbf{ and } $\sigma_{2} \neq 0$ \textbf{ then}\\
\hspace*{2cm} 5-Add the points with abscissa $x := \frac{-\sigma_{3}}{2 \sigma_{2}}$ to $S$.\\
 \hspace*{1cm} 6-end if\\
7-end for\\
8-Add the points at infinity  to $S$.\\
9- Return $S$.

\end{minipage}
}
\end{center}
\begin{Remark}
The cases 1 and 2 are not mutually  exclusive. Indeed, if a hyperelliptic curve \( C \) of even degree has at least one rational Weierstrass point, it is possible to find a hyperelliptic curve \( C' \) of odd degree  that is isomorphic to \( C \) over $\mathbb{Q}$. Thus, from the points on \( C' \), we can deduce the points on \( C \), see   Subsection \ref{sec:14}.
\end{Remark}
\section{Computation of isolated quadratic points on hyperelliptic curves of genus 2 and 3}\label{sec:51}

In this section, we extend the \texttt{Chabauty0} function \cite{Ca12, st}, as well as our function described in Section \ref{sec:4}, in order to compute the isolated quadratic points on hyperelliptic curves of genus 2 and 3.

\noindent Let $C$ be a hyperelliptic curve defined over $\mathbb{Q}$ whose Jacobian $J_C$ has rank 0. In this case, as discussed in the introduction, an isolated quadratic point is a point $P \in C(k)$ for some quadratic $k/\mathbb{Q}$  that does not lie in the preimage of $\mathbb{P}^1(\mathbb{Q})$ via a degree $2$ morphism $C \to \mathbb{P}^1$. The theorem of Bourdon, Ejder, Liu, Odumodu, and Viray \cite{BEV} asserts that such points are finite in number. It is therefore both natural and worthwhile to explicitly determine them.

\noindent   Suppose $P\in C(k)$ for some quadratic $k/\mathbb{Q}$ is isolated, let $G = \operatorname{Gal}(k/\mathbb{Q})$ denote its Galois group. We know that the divisor class $(P) +(\gamma(P)) - D_\infty$ lies in $J_C(\mathbb{Q})$, where $D_\infty$ denotes the divisor at infinity defined in Section \ref{sec:3}, and where $\sigma \in G$ denotes the Galois conjugation. Using our map described in Section \ref{sec:2}, the goal is to find the points in \( J_C(\mathbb{Q}) \) that are of the form
 \( \{P, \gamma
(P)\} \) for some quadratic $k/\mathbb{Q}$. Since we are working on the Kummer variety rather than directly on the Jacobian \( J_C \), it suffices to look for points on the Kummer variety that arise as images of divisors of the form \( \{P, \gamma
(P)\} \). We therefore obtain the following theorem

\begin{Theorem}\label{z2}
\begin{enumerate}
    \item Let \( R = (1 : \sigma_2 : \sigma_3 : \beta_0) \in K(\mathbb{Q}) \), where \( K \) denotes the Kummer surface associated with a genus 2 hyperelliptic curve. Then the preimage \( \kappa^{-1}(R) \) is of the form \( \{P, \gamma
(P)\} \) if and only if \( \sigma_3 \neq 0 \), and \( \sigma_2^2 - 4\sigma_3  \) is not a square in $\mathbb{Q}$. In the case, the \( x \)-coordinate of the point \( P \) is given by $$x = \frac{\sigma_2 + \sqrt{\sigma_2^2 - 4\sigma_3}}{2}$$.
    \item if \( R = (0 : \sigma_2 : \sigma_3 : \sigma_4 : \sigma_5 : \sigma_6 : \sigma_7 : \sigma_8) \in K(\mathbb{Q}) \), where \( K \) is the Kummer variety of a genus 3 hyperelliptic curve, then \( \kappa^{-1}(R) \) is of the form \( \{P, \gamma
(P)\} \), with $P$ is not quadratic Weierstrass point, if and only if \( \sigma_2 \neq 0 \), \( \sigma_4 \neq 0 \), and \( \sigma_3^2 - 4\sigma_2 \sigma_4 \) is not a square in $\mathbb{Q}$. The \( x \)-coordinate of the point \( P \) is given by $$x=\frac{-\sigma_3 + \sqrt{\sigma_3^2 - 4\sigma_2 \sigma_4}}{2\sigma_2}.$$
\end{enumerate}
\end{Theorem}
\begin{proof}
We provide a proof for case (2); case (1) can be proved in a similar manner, and in fact is simpler.

  (2)  Let \( C \) be a genus 3 hyperelliptic curve defined over \( \mathbb{Q} \) by the equation \( y^2 = f(x) \), where \( f \) is a polynomial of degree 7 or 8.  
Let \( k = \mathbb{Q}(\sqrt{d}) \) be a quadratic extension of \( \mathbb{Q} \), with \( d \in \mathbb{Q}^\times \) not a square in \( \mathbb{Q} \), and let \( G = \mathrm{Gal}(k/\mathbb{Q}) \) be its Galois group. Let 
$P := (a + b\sqrt{d}, y) \in C(k)$
with \( a \in \mathbb{Q} \), \( b \in \mathbb{Q}^\times \).\\ Suppose $Q\in K(\mathbb{Q})$ such that
$\kappa^{-1}(R) = \{P, \gamma
(P)\}$
where \( \gamma
 \in G \) is the nontrivial Galois automorphism.
Then,
$$ \kappa(Q) = (0 : 1 : -2a : a^2 - b^2 d : \ldots).$$ We aim to show that:
 \( \sigma_2 \ne 0 \),
 \( \sigma_4 \ne 0 \) and
 \( \sigma_3^2 - 4\sigma_2 \sigma_4 \notin (\mathbb{Q}^\times)^2\).
Assume for contradiction that \( \sigma_4 = 0 \). From the formula:
\[
\sigma_4 = x_1 x_2 = (a + b\sqrt{d})(a - b\sqrt{d}) = a^2 - b^2 d
\]
Thus, \( \sigma_4 = 0 \Rightarrow a^2 = b^2 d \Rightarrow a = \pm b \sqrt{d} \), which implies that:
either \( x_1 \) or \( x_2 \) is equal to 0. This contradicts the assumption that \( x_1 \notin \mathbb{Q} \). Hence, \( \sigma_4 \ne 0 \).
Now, observe that:
\[
\sigma_3 = - (x_1 + x_2) = -2a, \quad \sigma_2 = 1
\Rightarrow \sigma_3^2 - 4\sigma_2 \sigma_4 = 4a^2 - 4(a^2 - b^2 d) = 4b^2 d
\]

Since \( b \in \mathbb{Q}^\times \) and \( d \notin (\mathbb{Q}^\times)^2\), it follows that:
\[
4b^2 d \notin (\mathbb{Q}^\times)^2\Rightarrow \sigma_3^2 - 4\sigma_2 \sigma_4 \notin \mathbb{Q}^2
\]
\noindent  Let $R = (0 : \sigma_2 : \sigma_3 : \sigma_4 : \sigma_5 : \sigma_6 : \sigma_7 : \sigma_8) \in K(\mathbb{Q}).$
We associate to \( R \) the following quadratic polynomial:
\[
A(x) = \sigma_2 x^2 + \sigma_3 x + \sigma_4,
\]
with \( \sigma_2 \ne 0 \). Its discriminant is given by
\[
\Delta := \sigma_3^2 - 4\sigma_2 \sigma_4.
\]

Assume that \( \Delta \notin (\mathbb{Q}^\times)^2\). In particular, this implies \( \sigma_4 \ne 0 \), and the two roots of the equation are:
\[
x_1 = \frac{-\sigma_3 + \sqrt{\Delta}}{2\sigma_2}, \quad
x_2 = \frac{-\sigma_3 - \sqrt{\Delta}}{2\sigma_2}.
\]

These roots lie in the quadratic field \( k = \mathbb{Q}(\sqrt{\Delta}) \). Therefore, the rational divisor  associated to the pair \( (x_1, x_2) \) corresponds to two points \( P_1, P_2 \in C(k) \), where \( P_2 = \gamma
(P_1) \) for \( \gamma
 \in \mathrm{Gal}(k/\mathbb{Q}) \).

Hence, we conclude that:
\[
\kappa^{-1}(R) = \{P, \gamma
(P)\}.
\]

\end{proof}
\noindent The procedure for determining these points proceeds as follows: we begin by searching for quadratic Weierstrass points by searching for the quadratic irreducible factors of $f$. Then, within the set $TK \subseteq K(\mathbb{Q})$, we identify the points satisfying the conditions stated in the previous theorem.  Note that in this case, we do not obtain any points at infinity, since the $x$-coordinates of such points are, by definition, invariant under the action of the Galois group.


\section{Comparison with an  algorithm based on Coleman integration and generalizations }\label{sec:5}
\subsection{Comparison with the algorithm of Fernandez and Hashimoto }
Another approach one can  use to determine \( C(\mathbb{Q}) \), when \( C \) is a genus 3 hyperelliptic curve defined over \( \mathbb{Q} \) whose Jacobian has rank 0 and is irreducible, is the method of Fern\'andez and Hashimoto \cite{Fu12} implemented in Sage  \cite{degond9}. It works as follows:

Let \( J_C \) be the Jacobian of \( C \), and let \( p \) be a prime of good reduction for \( C \). A basis of the space of regular differentials \( H^0 (C_{\mathbb{Q}_p}, \Omega^1) \) is given by \( \{ \omega_0, \omega_1, \omega_2 \} \), where
$\omega_i = \frac{x^i}{2y} \, dx.$
For each \( i = 0,1,2 \), we define the  function $f_{i} : C(\mathbb{Q}_{p})\rightarrow \mathbb{Q}_{p}$ by :
\[
f_i (Q) = \int_{\infty}^{Q} \omega_i,
\]
 where the integral is a Coleman integral, see for instance \cite{MP2}.
For each point \( \tilde{Q} \) in \( C(\mathbb{F}_p) \), we compute the set of \( \mathbb{Q}_p \)-rational points \( Q \) reducing to \( \tilde{Q} \) such that \( f_i (Q) = 0 \) for \( i = 0,1,2 \).

If \( P \in J_C(\mathbb{Q}_p) \) is such that all integrals \( \int_P \omega_i = 0 \), then \( P \) belongs to \( J(\mathbb{Q}_p)_{\text{tors}} \) \cite[page 42]{sto}.
In particular, if \( Q \in C(\mathbb{Q}_p) \) is such that \( [Q - \infty] \) is torsion, then \( f_i(Q) = 0 \).
Thus, one only needs to find the points \( Q \) in the common zero set of the \( f_i \) and check whether \( Q \) is rational. It should be mentioned that these \( p \)-adic analytic computations can only be performed to finite precision, so a careful precision analysis is required. Here we discuss how this approach compares to ours.\\
The implementation of the algorithm by Fernandez and Hashimoto \cite{degond9} is limited to genus 3 hyperelliptic curves given by a degree 7 model. Indeed, for genus 3 hyperelliptic curves given by a degree 8 model with no rational Weierstrass points, there seems to be no available implementation of this algorithm. Thus, general even models could be handled by replacing the point at infinity with another rational point in the definition of \( f_i \) and using the \texttt{Magma}- implementation of Coleman integration \cite{BJen, BT1}. Our algorithm, on the other hand, is simpler and can be applied to any hyperelliptic curve model to compute the rational points as well as the isolated quadratic points. Since our code relies on the computation of Kummer varieties in genus 3 which in practice is   a notoriously complex task, it is generally slower than that of Fernández and Hashimoto. However, when the torsion subgroup is trivial, the use of Kummer varieties becomes unnecessary. In such cases, our method becomes significantly more efficient and clearly outperforms that of Fernández and Hashimoto. We discuss these execution times in detail in the examples that follow \ref{sec:14} and \ref{sec:15}. In  Subsection  \ref{sec:55},  we present an example of a hyperelliptic curve with no rational points, a case that cannot be handled by the algorithm \cite{Fu12} unless one combines it with an alternative approach such as the Mordell-Weil sieve \cite{Sieve}.
\subsection{ Possible generalizations of the algorithm}
\subsection*{ Hyperelliptic curves of genus $g\ge 4$}
The hyperelliptic curve used as input in the algorithm is assumed to have genus 3 and rank equal to 0. However, the rank computation becomes particularly complex when the genus of the curve is greater than or equal to 3. Nevertheless, if it is known with certainty that the rank of the Jacobian of our curve with genus \( g \geq 4 \) is 0, and if there is an explicit theory for the Kummer variety as well as a theory of heights for hyperelliptic curves of genus \( g \geq 4 \), it would then be possible to extend the algorithm of Müller and Reitsma for computing the rational torsion group of the Jacobian, \( J_C(\mathbb{Q})_{\text{tors}} \), and adapt our own algorithm. Such an explicit theory of the Kummer variety and of heights is work in progress due to Ludwig Fürst.


\subsection*{Genus 3 non-hyperelliptic curves }
The algorithm we have developed in this paper requires the existence of an explicit theory of Kummer varieties and heights for algebraic curves. Unfortunately, for the case of non-hyperelliptic curves of genus 3, such an explicit theory does not exist. The work of van Bommel in \cite{van,van1} allows for the computation of the rational torsion subgroup in the case of non-hyperelliptic curves of genus 3. But, the problem of rank computation in this case is complex, even though an algorithm by Bruin, Poonen, and Stoll, described in \cite{BPS}, exists. However, no general implementation is available. In general, if we know with certainty that a non-hyperelliptic curve of genus 3 has an irreducible Jacobian and rank 0, we can use van Bommel's algorithm to compute \( J_C(\mathbb{Q})_{\text{tors}} \), and then apply the strategy described in Section \ref{sec:2}.
If the Jacobian is not irreducible, this means that \( C \) is a covering of a curve \( E \) of genus 1, whose rank is also 0. In this case, it is preferable to apply the rank 0 quotient strategy, as described in \cite{EBR}.
\subsection*{Points of degree \( d > 2 \) on a hyperelliptic curve of genus 3}
Let \( C \) be a hyperelliptic curve of genus \( 3 \), defined over \( \mathbb{Q} \), and assume it admits at least one rational Weierstrass point. In this case, \( C \) has a unique point at infinity.
\noindent Let \( k/\mathbb{Q} \) be a Galois extension of degree \( d \geq 3 \), and let \( G = \mathrm{Gal}(k/\mathbb{Q}) \) be its Galois group. We can extend our algorithm to compute isolated points of degree \( d \) by constructing the class of a divisor \( Q \in J_C(\mathbb{Q}) \) Let 
\(
Q = [D - d\infty],
\)
where \(D\) is an effective divisor of degree \(d\). 
In the case where \(G\) is cyclic, we retain this degree \(d\); 
otherwise, we consider a divisor of degree \(n\), 
where \(n\) is a strictly positive integer dividing \(d\). We are currently working on the case \(d = 3\).
\section{Data and Examples}\label{sec:6}
We ran the \texttt{Magma} implementation \cite{degond8} of our algorithm on a database of 7396 genus 3 hyperelliptic curves whose Jacobian has rank 0, taken from the list \cite{MB12}. This list \cite{MB12} corresponds to the set of all rank 0 curves in Sutherland's database of genus 3 hyperelliptic curves with discriminant bounded by $10^7$ \cite{degond7}.

The rank computation was performed using \texttt{Magma}’s \textbf{\texttt{RankBound}} function, which provides an upper bound for the rank of the curve. For all these curves, the upper bound returned by \textbf{\texttt{RankBound}}  was always 0. Among these 7396 hyperelliptic curves, 4342 are given by an odd-degree model and 3054 are given by an even-degree model, among which 1576 have at least one Weierstrass point.
Our implementation proves that, for each of the studied curves, the entire set of rational points
equals the set of rational points of naive height less than \(10^5\).

\begin{figure}[h!]
\centering
 \includegraphics[width=0.7\textwidth]{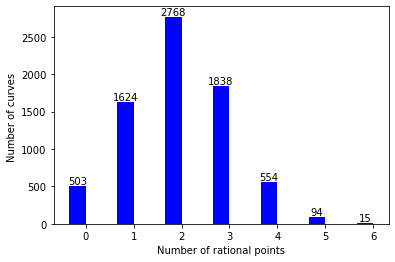}

 \caption{ Graph of number of rational points }
 \label{fig:f1}
\end{figure}

\begin{figure}[h!]
\centering
 \includegraphics[width=0.7\textwidth]{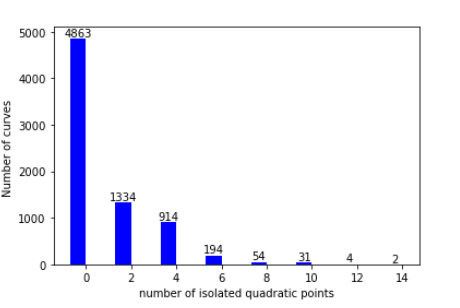}

 \caption{ Graph of number of isolated quadratic points in genus 3 }
 \label{fig:f3}
\end{figure}

\begin{figure}[h!]
\centering
 \includegraphics[width=0.7\textwidth]{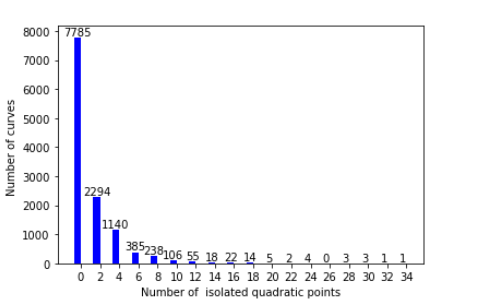}

 \caption{ Graph of number of isolated quadratic points in genus 2 }
 \label{fig:f2}
\end{figure}

We show in Figure~\ref{fig:f1} below how many of the curves in our database have a certain number of
rational points. We observe that the maximum number of points is six, and that a vast majority of the curves have three or fewer rational points.
 These full sets of rational points  can be found in the master list of Chabauty0Genus3 data on the 7396 hyperelliptic curves \textbf{allpointsdata.txt} available at \href{https://github.com/Brice202145/Chabauty0Genus3}{https://github.com/Brice202145/Chabauty0Genus3}.
Furthermore, we also investigated whether the curves in our database contain any isolated quadratic points. By running our code, we obtained the results summarized in the figure ~\ref{fig:f3}  .
We queried the LMFDB database to retrieve hyperelliptic curves of genus~2 whose Jacobians has rank~0. We tested our algorithm on 12,075 such genus~2 curves. For more details, see the following link: \url{https://www.lmfdb.org/Genus2Curve/Q/?analytic_rank=0}. We used our implementation to compute the set \( \mathcal{ L} \) of isolated quadratic points. Our implementation, applied to these 12,075 curves, yields the results shown in  Figure ~\ref{fig:f2}.

We present here some computation times illustrating the difference between our algorithm and the one described in \cite{Fu12}. The computations in this subsection were carried out on a computer with an Intel\textsuperscript{\textregistered} Core\textsuperscript{TM} i5-7300U CPU @ 2.60\,GHz × 4 and 8.0\,GiB of RAM
\subsection{An example from \cite{Fu12}}\label{sec:14}

Let \( C/\mathbb{Q} \) be the hyperelliptic curve of genus 3 defined by the equation:
\[C \colon y^2 = x^8 + 4x^7 + 6x^6 + 4x^5 - 7x^4 -16x^3 + 8x.\]
Note that this model is a simplified model of the equation given in \cite[Example 3]{Fu12}. In \cite{Fu12}, the authors show that the number of rational points on this curve  equals Stoll's bound in $p:=7$. Indeed, when \( p > 2g \) and the curve has good reduction modulo \( p \), and if \( r < g-1 \), Stoll \cite[Corollary 6.7]{Stoll2} has shown that:
$\#C(\mathbb{Q}) < \#C(\mathbb{F}_p) + 2r.$
In our case, \( r = 0 \), so we have:
$\#C(\mathbb{Q}) < \#C(\mathbb{F}_p).$
 We will now recompute \( C(\mathbb{Q}) \) using our algorithm.\\
The curve \( C \) admits at least one rational Weierstrass point, allowing us to find an odd-degree model defined over \( \mathbb{Q} \). Such a model is given by:
\[
C' \colon y^2 = 8x^7 - 16x^5 - 7x^4 + 4x^3 + 6x^2 + 4x + 1,
\]
which admits three rational Weierstrass points: \((-1, 0) \), \((-1/2, 0)\), and \((1, 0)\). The set $TK$ is given by: \\
$ TK=\{ (1/2 : 0 : 0 : 0 : -4 : 0 : 1 : 0), (0 : 1/2 : 0 : 0 : 0 : 0 : 0 : 1), (0 : 4
: 6 : 2 : 7 : 3 : 1 : 0), (0 : 0 : 0 : 0 : -1/8 : 1/8 : -1/8 : 1), (0 : 1/6 : 0
: -1/6 : 1/6 : 0 : 1/6 : 1), (-1/8 : 1/2 : 3/2 : 1/2 : 3/2 : 1/2 : -1/4 : 1), (0
: 4 : 2 : 0 : 1 : 0 : 0 : 0), (-1/32 : 1/2 : 0 : -1/4 : 0 : -1/4 : -1/16 : 1),
(-1/8 : -3/2 : -1/2 : 0 : 1 : 1/2 : 1/4 : 1), (1/24 : 1/6 : -1/3 : -1/6 : -1/6 :
0 : 1/12 : 1), (-1/32 : 0 : 1/4 : 0 : 1/4 : 1/8 : 1/16 : 1), (1/64 : 0 : -1/4 :
-1/8 : -1/8 : 0 : 1/32 : 1), (-1/8 : 1/2 : 1/2 : 0 : 2 : 1/2 : 1/4 : 1), (-1/16
: 1/2 : 0 : 0 : 3/2 : 1/2 : -1/8 : 1), (0 : 0 : 0 : 0 : 1 : 0 : 0 : 0), (-1/16 :
1/2 : 1/2 : 0 : 1 : 1/4 : 1/8 : 1), (0 : -1/6 : 1/6 : 0 : -1/6 : 0 : 0 : 1),
(1/48 : -1/6 : -1/3 : -1/6 : -1/6 : 0 : 1/24 : 1), (-1/16 : 1/2 : 1/2 : 0 : 1/2
: 0 : -1/8 : 1), (0 : 0 : 0 : 0 : 1 : 1/2 : 1/4 : 1), (1/16 : 1/2 : 0 : 0 : 0 :
0 : 1/8 : 1), (1/8 : 1/2 : 0 : 0 : 0 : 1/2 : 1/4 : 1), (0 : 0 : 0 : 0 : 1/8 :
1/8 : 1/8 : 1), (0 : 1/4 : -1/8 : -1/8 : 3/16 : 1/16 : 1/16 : 1), (0 : 1/10 :
-1/10 : -1/10 : 1/5 : 1/10 : 1/10 : 1), (-1/64 : 0 : -1/8 : -1/8 : 1/4 : 1/8 :
3/32 : 1), (0 : 1/2 : 1/2 : 0 : 1/2 : 0 : 0 : 1) \}$

We now search for \( R \in \text{TK} \) that correspond to rational affine non-Weierstrass points on  \( C' \). We find exactly one such \( R \):
\[
R = (0 : \frac{1}{2} : 0 : 0 : 0 : 0 : 0 : 1),
\]
which corresponds to the rational points \((0, 1)\) and \((0, -1 )\) on \( C' \).\\
Hence,
\[ C'(\mathbb{Q})=\{(0, -1 ), (0, 1 ),(-1, 0) , (-1/2, 0), (1, 0), \infty \},  \] and so
\[  C(\mathbb{Q}) =\{ (-2, 0), (-1, 0 ), (0, 0 ), (1, 0), \infty_{-},\infty_{+} \}.
\]
Taking the reduction modulo 7, we also find 6 points over \( \mathbb{F}_7 \):
\[
 C(\mathbb{F}_7)=\{(\infty_{1}, \infty_{2}, (0, 0), (1, 0), (5, 0 ), (6, 0) \}.
\]
By Stoll’s theorem, we know that there cannot be any more rational points. Thus, this bound is sharp for \( p = 7 \).

\noindent The computation using our function \texttt{Chabauty0Genus3} for the curve \( C' \) took 8.5 seconds, whereas the use of the \texttt{chabauty\_coleman} function required only 0.5 seconds. This difference is explained by the fact that the torsion subgroup of the Jacobian of \( C' \) is large, with a torsion bound equal to 48  and that the Kummer variety is a complicated object, so that computing with it is expensive.

By examining the set \( TK \), we detected the point \( R \in TK \) that corresponds to isolated quadratic points on the curve \( C' \). This point is given by:
\[
R := (0 : \tfrac{1}{10} : -\tfrac{1}{10} : -\tfrac{1}{10} : \tfrac{1}{5} : \tfrac{1}{10} : \tfrac{1}{10} : 1).
\]
It satisfies the conditions of the theorem~\ref{z2}, and our implementation returns the set:\\
\[ 
\left\{
\begin{aligned}
&\left( \tfrac{1}{2}(-\sqrt{5} + 1), \tfrac{1}{2}(-3\sqrt{5} + 7) \right),
&&\left( \tfrac{1}{2}(-\sqrt{5} + 1), \tfrac{1}{2}(3\sqrt{5} - 7) \right), \\
&\left( \tfrac{1}{2}(\sqrt{5} + 1), \tfrac{1}{2}(-3\sqrt{5} - 7) \right),
&&\left( \tfrac{1}{2}(\sqrt{5} + 1), \tfrac{1}{2}(3\sqrt{5} + 7) \right)
\end{aligned}
\right\}
\]
\subsection{Genus 3 hyperelliptic curve of odd degree}\label{sec:15}
Let $C$ be the curve of genus $g=3$ defined over $\mathbb{Q}$ by the equation $$C \colon y^2=x^7 - 47x^6 - 86x^5 + 49x^4 + 94x^3 - 90x^2 + 28x - 3.$$
Let $J_C$ be the Jacobian of $C$. The  \textbf{\texttt{RankBound}} function in  \texttt{Magma}  returns $0$ so the Mordell Weil rank of $J_C$ is $0$ over $\mathbb{Q}$. Then, by the Mordell-Weil theorem, $J_C(\mathbb{Q})=J_C(\mathbb{Q})_\text{tors}$
where $ J_C(\mathbb{Q})_\text{tors}$ is the rational torsion subgroup of $J_C$. Calculating $J_C(\mathbb{F}_p)$ for a few primes of good reduction \([3, 5, 7, 11, 17, 23, 31, 41, 43]\) shows that $J_C(\mathbb{Q})_{\mathrm{tors}}$ is trivial. Hence, we find that $TK$ is trivial.  Applying our code, the set of rational points on $C$ is
$$C(\mathbb{Q})=\{\infty\}.$$
Our \texttt{Chabauty0Genus3} function finished in 0.1 seconds, whereas the \texttt{Chabauty\_Coleman} function took 2 seconds.

\subsection{Genus 3 hyperelliptic curve of even degree without rational  Weierstrass points}\label{sec:16}
Consider the hyperelliptic curve \( C \) defined over \( \mathbb{Q} \) by the equation:
\[
C \colon y^2 =(x^2 - x + 1)(x^6 + x^5 - 6x^4 - 3x^3 + 14x^2 - 7x + 1).
\]
Let \( J \) denote its Jacobian, and let \( K \) be the Kummer variety of \( J \).

Let \( p = 5 \) be a prime of good reduction for both \( C \) and \( J \). Using \texttt{\texttt{Magma}}, we compute \( J(\mathbb{F}_5) \). For each point \( P \in J(\mathbb{F}_5) \), we apply our code, which computes the set \( TK \) which contains 50 points.\\
 Our code identifies the rational points \( R \in TK \) that lift to divisor classes of the form \( \{P, P\} \) in \( J(\mathbb{Q})_{\text{tors}} \), where \( P \in C(\mathbb{Q}) \). We have found the following points
 \[
R = (0 : 1/6 : 0 : 0 : 0 : 0 : 0 : 1) \quad \text{and} \quad R = (0 : 1/4 : -1/2 : 1/4 : 3/4 : -1/2 : 1/4 : 1),
\]
correspond to the rational points \( (1, -1) \), \( (1, -1) \), \( (0, 1) \), and \( (0, -1) \) on the curve \( C \).

Including the two rational points at infinity, the full set of rational points on \( C \) is given by:
\[
C(\mathbb{Q}) = \{ \infty_1, \infty_2, (-1, 1), (1, -1), (0, -1), (0, 1) \}.
\]

Our code also detects the points \( R \in TK \) that correspond to divisor classes of the form \( \{P, \sigma(P)\} \) in \( J(\mathbb{Q})_{\text{tors}} \), where \( P \) is an isolated quadratic point. These are exactly the following:
\[
\begin{aligned}
&(0 : 2/9 : -4/9 : 1/9 : 7/9 : -2/9 : 1/18 : 1), \\
&(0 : 1/4 : 1/2 : -1/4 : 5/4 : -1/2 : 1/4 : 1), \\
&(0 : -1/18 : 0 : 1/9 : -1/9 : 0 : -2/9 : 1),
\end{aligned}
\]
which correspond to the following set of  isolated quadratic points over \( \mathbb{Q}(\sqrt{2}) \):
\[ 
\left\{
\begin{aligned}
&\left( \tfrac{1}{2}(-\sqrt{2} + 2), \tfrac{1}{4}(4\sqrt{2} - 5) \right),
&&\left( \tfrac{1}{2}(-\sqrt{2} + 2), \tfrac{1}{4}(-4\sqrt{2} + 5) \right),
&&\left( \tfrac{1}{2}(\sqrt{2} + 2), \tfrac{1}{4}(-4\sqrt{2} - 5) \right), \\
&\left( \tfrac{1}{2}(\sqrt{2} + 2), \tfrac{1}{4}(4\sqrt{2} + 5) \right),
&&(\sqrt{2} - 1, 8\sqrt{2} - 11),
&&(\sqrt{2} - 1, -8\sqrt{2} + 11), \\
&(-\sqrt{2} - 1, 8\sqrt{2} + 11),
&&(-\sqrt{2} - 1, -8\sqrt{2} - 11),
&&(-\sqrt{2}, -4\sqrt{2} - 5), \\
&(-\sqrt{2}, 4\sqrt{2} + 5),
&&(\sqrt{2}, -4\sqrt{2} + 5),
&&(\sqrt{2}, 4\sqrt{2} - 5)
\end{aligned}
\right\}
\]
Finally, our code also detects the Weierstrass isolated quadratic points defined over \( \mathbb{Q}(\sqrt{-3}) \), which are:
\[
\left( \tfrac{1}{2}(\sqrt{-3} + 1), 0 \right) \quad \text{and} \quad \left( \tfrac{1}{2}(-\sqrt{-3} + 1), 0 \right).
\]
\subsection{Even-degree hyperelliptic curve without $\mathbb{Q}$-points at infinity}\label{54}

Let \( C/\mathbb{Q} \) be a genus 3 hyperelliptic curve given by the equation:
\[
C\colon y^2 = f(x) = 5x^8 - 22x^7 + 53x^6 - 74x^5 + 52x^4 + 2x^3 - 11x^2 + 2x + 1.
\]
Since the leading coefficient of \( f(x) \) is not a square in \( \mathbb{Q} \), the curve \( C \) has no rational points at infinity. We then seek the affine rational points on \( C \).

To do so, we compute the set \( TK \) of points in the Kummer variety \( K(\mathbb{Q}) \) that lift to \( J_C(\mathbb{Q}) \), where \( J_C \) is the Jacobian of \( C \). Let \( p := 5 \), a prime of good reduction for \( J_C \). After running our code for all \( P \in J_C(\mathbb{F}_p) \), we find that the set $TK$  is equal to:\\
$TK = \{
 (-1/448 : -5/64 : 1/448 : 1/112 : 9/448 : -3/224 : 1/32 : 1),(0 : -1/12 : 0 :0 : 0 : 0 : 0 : 1),
(-1/156 : -1/39 : -5/78 : -1/78 : -4/39 : -7/78 : 1/13 : 1)\}.$

The only point \( R \in TK \) satisfying \( \sigma_1 = 0 \), \( \sigma_2 \neq 0 \), and \( \sigma_3^2 - 4\sigma_2\sigma_4 = 0 \) is:
\[
R = \left( 0 : -\frac{1}{12} : 0 : 0 : 0 : 0 : 0 : 1 \right).
\]
Thus, the affine rational points of \( C(\mathbb{Q}) \) that are not Weierstrass points have abscissa \( x = 0 \), which corresponds to the points \( (0,1) \) and \( (0,-1) \).

The curve \( C \) has no rational Weierstrass points, hence the set of rational points is given by:
\[
C(\mathbb{Q}) = \{ (0,1), (0,-1) \}.
\]
\subsection{A hyperelliptic curve without rational points}\label{sec:55}
Let \( C \) be the hyperelliptic curve defined by the equation:
\[
C \colon y^2 = f(x) = -4x^8 + 8x^7 - 3x^6 - 16x^5 + 26x^4 - 16x^3- 3x^2 + 8x - 4.
\]
We observe that the polynomial \( f(x) \) has no roots in \( \mathbb{Q} \). Therefore, the curve \( C \) does not admit any rational Weierstrass points.
Moreover, the leading coefficient of \( f(x) \) is \(-4\), which is not a square in \( \mathbb{Q} \). As a result, the curve \( C \) also has no rational points at infinity.
We then focus on affine rational points \( P \in C(\mathbb{Q}) \) that are not Weierstrass points. To this end, we compute the set
\[
TK = \{
( -1/12 : 1 : -2/3 : 2/3 : -5/6 : -2/3 : 1 : 0 ),
( 0 : 1/2 : -1/2 : 1/2 : 0 : -1/2 : 1/2 : 1 )
\}.
\]

The point of interest is:
\[
R = \left( 0 : 1/2 : -1/2 : 1/2 : 0 : -1/2 : 1/2 : 1 \right).
\]

However, the coordinates of this point do not satisfy the condition:
\[
\sigma^{2}_{3} = 4\sigma_2 \sigma_4.
\]

Therefore, the curve \( C \) does not have any  no rational points:
\[
C(\mathbb{Q}) = \varnothing.
\]
The coordinates of the point \( R \) satisfy the conditions of the theorem \ref{z2}. Therefore, it corresponds to the following affine isolated quadratic  points :
\[
\left( \tfrac{1}{2}(-\sqrt{-3} + 1), \tfrac{1}{2}(\sqrt{-3} + 1) \right), \quad
\left( \tfrac{1}{2}(-\sqrt{-3} + 1), \tfrac{1}{2}(-\sqrt{-3} - 1) \right),
\]
\[
\left( \tfrac{1}{2}(\sqrt{-3} + 1), \tfrac{1}{2}(\sqrt{-3} - 1) \right), \quad
\left( \tfrac{1}{2}(\sqrt{-3} + 1), \tfrac{1}{2}(-\sqrt{-3} + 1) \right).
\]

We then investigate whether there are any quadratic Weierstrass isolated points. Our code returns the following two points:
\[
\left( \tfrac{1}{2}(2\sqrt{-7} + 3), 0 \right), \quad
\left( \tfrac{1}{2}(-2\sqrt{-7} + 3), 0 \right).
\]
\subsection{ A genus 2  curve that  does not admit quadratic  Weierstrass points}

Let $C$ be the hyperelliptic curve of genus 2 defined over $\mathbb{Q}$ by the equation
\[ y^2 = 8x^5 - 7x^4 + 6x^3 - x^2 + 2x + 1. \]
The function \texttt{\texttt{RankBound}} returns that the rank of $J_C$ over $\mathbb{Q}$ is 0.

We use our code to compute the set $\mathcal{L}$ of quadratic points. In the first step, our code searches for quadratic  Weierstrass points by finding irreducible polynomials $g$ over $\mathbb{Q}$ that divide the polynomial $8x^5 - 7x^4 + 6x^3 - x^2 + 2x + 1$. In this case, such a polynomial does not exist because the polynomial is irreducible over $\mathbb{Q}$.

In the second step, we compute the  isolated quadratic  affine points. To do this, we compute the set $TK$, where $K$ is the Kummer surface of the Jacobian $J_C$. It is given by
\begin{align*}
TK = \{ & (1 : 2 : 1 : -6),\ (0 : 1 : 1 : 8),\ (1 : 1 : 0 : -2),\ (2 : 0 : 1 : 4),\ (0 : 1 : 0 : 0),\\
        & (1 : -1 : 1 : 18),\ (1 : 1 : 0 : 10),\ (1 : 3 : -1 : 2),\ (1 : 1 : 1 : 2),\ (1 : 0 : 0 : 2) \}.
\end{align*}

The points $R$ in $TK$ that we are interested in are:
\[(1 : 2 : 1 : -6),\ (2 : 0 : 1 : 4),\ (1 : -1 : 1 : 18),\ (1 : 3 : -1 : 2), \text{and} \ (1 : 1 : 1 : 2),\]
which correspond to the isolated quadratic points on $C$ as follows:
\begin{align*}
\mathcal{L} = \{ & (1/2\sqrt{-2},\ \tfrac{1}{2}(-\sqrt{-2} - 1)),\ (-1/2\sqrt{-2},\ \tfrac{1}{2}(-\sqrt{-2} + 1)),\\
                 & (1/2\sqrt{-2},\ \tfrac{1}{2}(\sqrt{-2} + 1)),\ (-1/2\sqrt{-2},\ \tfrac{1}{2}(\sqrt{-2} - 1)),\\
                 & (\tfrac{1}{2}(-\sqrt{-3} - 1),\ -\sqrt{-3} - 3),\ (\tfrac{1}{2}(-\sqrt{-3} - 1),\ \sqrt{-3} + 3),\\
                 & (\tfrac{1}{2}(\sqrt{-3} - 1),\ -\sqrt{-3} + 3),\ (\tfrac{1}{2}(\sqrt{-3} - 1),\ \sqrt{-3} - 3),\\
                 & (\tfrac{1}{2}(-\sqrt{-3} + 1),\ 2),\ (\tfrac{1}{2}(-\sqrt{-3} + 1),\ -2),\\
                 & (\tfrac{1}{2}(\sqrt{-3} + 1),\ 2),\ (\tfrac{1}{2}(\sqrt{-3} + 1),\ -2),\\
                 & (\tfrac{1}{2}(-\sqrt{13} + 3),\ -7\sqrt{13} + 25),\ (\tfrac{1}{2}(-\sqrt{13} + 3),\ 7\sqrt{13} - 25),\\
                 & (\tfrac{1}{2}(\sqrt{13} + 3),\ 7\sqrt{13} + 25),\ (\tfrac{1}{2}(\sqrt{13} + 3),\ -7\sqrt{13} - 25) \}.
\end{align*}

\subsection{ A genus 2 curve that admits quadratic  Weierstrass points}

Let $C$ be the hyperelliptic curve of genus 2 defined over $\mathbb{Q}$ by the equation 
\[ y^2 = f(x) = x^6 - 22x^4 + 157x^2 - 360. \]
The polynomial $f(x)$ factors over $\mathbb{Q}$ as 
\[ f(x) = (x - 3)(x + 3)(x^2 - 8)(x^2 - 5), \]
which yields the following quadratic  Weierstrass points:
\[(2\sqrt{2}, 0), (-2\sqrt{2}, 0), (-\sqrt{5}, 0), (\sqrt{5}, 0).\]

We now search for non-Weierstrass quadratic isolead points by computing the set $TK$:
\[
\begin{aligned}
TK = \{ &(0 : 0 : 1 : 22), (1 : 0 : -7 : -206), (1 : 0 : -6 : -192), (0 : 1 : -3 : -54), \\
&(2 : 1 : -13 : -394), (2 : -1 : -13 : -394), (0 : 1 : 3 : 54), (1 : 0 : -10 : -256), \\
&(1 : 0 : -8 : -221), (1 : 0 : -9 : -238), (1 : 0 : -11 : -278), (1 : 0 : -5 : -182), \\
&(2 : 0 : -15 : -426) \}.
\end{aligned}
\]

We do not find any $R \in TK$ corresponding to rational points. However, except for those $R$ whose first coordinate is zero, the remaining elements correspond to quadratic isolead points on $C$.

Thus, the set $\mathcal{L}$ of all such quadratic isolead points is:
\[
\begin{aligned}
\mathcal{L} = \{ 
&(1/4(-\sqrt{105} + 1), 1/16(\sqrt{105} + 19)), (1/4(-\sqrt{105} + 1), 1/16(-\sqrt{105} - 19)), \\
&(1/4(\sqrt{105} + 1), 1/16(\sqrt{105} - 19)), (1/4(\sqrt{105} + 1), 1/16(-\sqrt{105} + 19)), \\
&(1/4(-\sqrt{105} - 1), 1/16(\sqrt{105} - 19)), (1/4(-\sqrt{105} - 1), 1/16(-\sqrt{105} + 19)), \\
&(1/4(\sqrt{105} - 1), 1/16(\sqrt{105} + 19)), (1/4(\sqrt{105} - 1), 1/16(-\sqrt{105} - 19)), \\
&(-\sqrt{7}, 2), (-\sqrt{7}, -2), (\sqrt{7}, 2), (\sqrt{7}, -2), \\
&(-\sqrt{11}, 6), (-\sqrt{11}, -6), (\sqrt{11}, 6), (\sqrt{11}, -6), \\
&(\sqrt{10}, -\sqrt{10}), (\sqrt{10}, \sqrt{10}), (-\sqrt{10}, -\sqrt{10}), (-\sqrt{10}, \sqrt{10}), \\
&(\sqrt{6}, -\sqrt{6}), (\sqrt{6}, \sqrt{6}), (-\sqrt{6}, -\sqrt{6}), (-\sqrt{6}, \sqrt{6}), \\
&(1/2\sqrt{30}, -1/4\sqrt{30}), (1/2\sqrt{30}, 1/4\sqrt{30}), \\
&(-1/2\sqrt{30}, -1/4\sqrt{30}), (-1/2\sqrt{30}, 1/4\sqrt{30}), \\
&(2\sqrt{2}, 0), (-2\sqrt{2}, 0), (-\sqrt{5}, 0), (\sqrt{5}, 0)
\}.
\end{aligned}
\]


\begin{thebibliography}{10}

\bibitem{BBC}
Jennifer~S. Balakrishnan, Francesca Bianchi, Victoria Cantoral-Farf{\'a}n,
  Mirela {\c{C}}iperiani, and Anastassia Etropolski.
\newblock Chabauty-{Coleman} experiments for genus 3 hyperelliptic curves.
\newblock In {\em Research directions in number theory. Women in numbers IV.
  Proceedings of the women in numbers, WIN4 workshop. Banff International
  Research Station, Banff, Alberta, Canada, August 14--18, 2017}, pages 67--90.
  Cham: Springer, 2019.

\bibitem{Bjen12}
Jennifer~S. Balakrishnan, Robert~W. Bradshaw, and Kiran~S. Kedlaya.
\newblock Explicit {Coleman} integration for hyperelliptic curves.
\newblock In {\em Algorithmic number theory. 9th international symposium,
  ANTS-IX, Nancy, France, July 19--23, 2010. Proceedings}, pages 16--31.
  Berlin: Springer, 2010.

\bibitem{BT1}
Jennifer~S. Balakrishnan and Jan Tuitman.
\newblock Magma code.
\newblock Available at \url{ https://github.com/jtuitman/Coleman}.

\bibitem{BJen}
Jennifer~S. Balakrishnan and Jan Tuitman.
\newblock Explicit {Coleman} integration for curves.
\newblock {\em Math. Comput.}, 89(326):2965--2984, 2020.

\bibitem{van}
Raymond~van Bommel.
\newblock Computing torsion for plane quartics without using height bounds.
\newblock Preprint, {arXiv}:2301.08312 [math.{NT}], 2023.

\bibitem{van1}
Raymond~van Bommel.
\newblock genus3torsion. magma code, 2025.
\newblock \url{https://github.com/rbommel/genus3torsion}.

\bibitem{BEV}
Abbey Bourdon, {\"O}zlem Ejder, Yuan Liu, Frances Odumodu, and Bianca Viray.
\newblock On the level of modular curves that give rise to isolated
  {{\(j\)}}-invariants.
\newblock {\em Adv. Math.}, 357:33, 2019.
\newblock Id/No 106824.

\bibitem{BPS}
Nils Bruin, Bjorn Poonen, and Michael Stoll.
\newblock Generalized explicit descent and its application to curves of genus
  3.
\newblock {\em Forum Math. Sigma}, 4:80, 2016.
\newblock Id/No e6.

\bibitem{Sieve}
Nils Bruin and Michael Stoll.
\newblock The {Mordell}-{Weil} sieve: proving non-existence of rational points
  on curves.
\newblock {\em LMS J. Comput. Math.}, 13:272--306, 2010.

\bibitem{Ca12}
John W.~S. Cassels and Victor~E. Flynn.
\newblock {\em Prolegomena to a middlebrow arithmetic of curves of genus 2},
  volume 230 of {\em Lond. Math. Soc. Lect. Note Ser.}
\newblock Cambridge: Cambridge Univ. Press, 1996.

\bibitem{cl}
Claude Chabauty.
\newblock Sur les points rationnels des courbes alg{\'e}briques de genre
  sup{\'e}rieur {\`a} l'unit{\'e}.
\newblock {\em C. R. Acad. Sci., Paris}, 212:882--885, 1941.

\bibitem{CR}
Robert~F. Coleman.
\newblock Effective {Chabauty}.
\newblock {\em Duke Math. J.}, 52:765--770, 1985.

\bibitem{degond9}
Mar{\'{\i}}a~I. de~Frutos-Fern{\'a}ndez and Sachi Hashimoto.
\newblock Sage code.
\newblock
  \url{https://github.com/sachihashimoto/rational-points-hyperelliptic.}

\bibitem{Fu12}
Mar{\'{\i}}a~I. de~Frutos-Fern{\'a}ndez and Sachi Hashimoto.
\newblock Computing rational points on rank 0 genus 3 hyperelliptic curves.
\newblock In {\em Arithmetic geometry, number theory, and computation}, pages
  449--460. Cham: Springer, 2021.

\bibitem{EBR}
Tony Ezome, Brice~Miayoka Moussolo, and Régis~Babindamana Freguin.
\newblock Rational {Points} of some genus $3$ curves from the rank $0$ quotient
  strategy.
\newblock {\em INTEGERS}, 25:2, 2025.

\bibitem{FG}
Gerd Faltings.
\newblock Finiteness theorems for abelian varieties over number fields.
\newblock {\em Invent. Math.}, 73:349--366, 1983.

\bibitem{zbMATH01466163}
Marc Hindry and Joseph~H. Silverman.
\newblock {\em Diophantine geometry. {An} introduction}, volume 201 of {\em
  Grad. Texts Math.}
\newblock New York, NY: Springer, 2000.

\bibitem{Kla9}
Matthew~James Klassen.
\newblock {\em Algebraic points of low degree on curves of low rank}.
\newblock Ph.d. thesis, The University of Arizona, Ann Arbor, MI, 1993.

\bibitem{LLL}
A.~K. Lenstra, H.~W. Lenstra, and L{\'a}szl{\'o} Lov{\'a}sz.
\newblock Factoring polynomials with rational coefficients.
\newblock {\em Math. Ann.}, 261:515--534, 1982.

\bibitem{MP2}
William McCallum and Bjorn Poonen.
\newblock The method of {Chabauty} and {Coleman}.
\newblock In {\em Explicit methods in number theory. Rational points and
  Diophantine equations.}, pages 99--117. Paris: Soci{\'e}t{\'e}
  Math{\'e}matique de France (SMF), 2012.

\bibitem{degond8}
Brice Miayoka.
\newblock Chabauty0 on hyperelliptic curve of genus 3.
\newblock Magma code \url{https://github.com/Brice202145/Chabauty0Genus3}.

\bibitem{MB12}
Brice Miayoka.
\newblock Genus 3 hyperelliptic curves with {M}ordell-{W}eil rank 0. 2025,.
\newblock see \url{https://github.com/Brice202145/genus3-hyperellipticcurves}.

\bibitem{Mur12}
J.~Steffen M{\"u}ller.
\newblock Explicit {Kummer} varieties of hyperelliptic {Jacobian} threefolds.
\newblock {\em LMS J. Comput. Math.}, 17:496--508, 2014.

\bibitem{Mur13}
J.~Steffen M{\"u}ller and Berno Reitsma.
\newblock Computing torsion subgroups of {Jacobians} of hyperelliptic curves of
  genus 3.
\newblock {\em Res. Number Theory}, 9(2):26, 2023.
\newblock Id/No 23.

\bibitem{zbMATH01714029}
Bjorn Poonen.
\newblock Computing torsion points on curves.
\newblock {\em Exp. Math.}, 10(3):449--465, 2001.

\bibitem{degond6}
Berno Reitsma.
\newblock Computing torsion subgroups of jacobians of hyperelliptic curves of
  genus 3 over $\mathbb{Q}$.
\newblock see \url{https://github.com/bernoreitsma/g3hyptorsion}.

\bibitem{Sisk}
Samir Siksek.
\newblock Chabauty for symmetric powers of curves.
\newblock {\em Algebra Number Theory}, 3(2):209--236, 2009.

\bibitem{st}
Michael Stoll.
\newblock Magma-related directory.
\newblock See \url{http://www.mathe2.uni-bayreuth.de/stoll/magma/index.html}.

\bibitem{Rk}
Michael Stoll.
\newblock Implementing 2-descent for {Jacobians} of hyperelliptic curves.
\newblock {\em Acta Arith.}, 98(3):245--277, 2001.

\bibitem{Stoll}
Michael Stoll.
\newblock On the height constant for curves of genus two. {II}.
\newblock {\em Acta Arith.}, 104(2):165--182, 2002.

\bibitem{Stoll2}
Michael Stoll.
\newblock Independence of rational points on twists of a given curve.
\newblock {\em Compos. Math.}, 142(5):1201--1214, 2006.

\bibitem{st12}
Michael Stoll.
\newblock An explicit theory of heights for hyperelliptic {Jacobians} of genus
  three.
\newblock In {\em Algorithmic and experimental methods in algebra, geometry,
  and number theory}, pages 665--715. Cham: Springer, 2017.

\bibitem{sto}
Michael Stoll.
\newblock {\em Arithmetic of Hyperelliptic Curves}.
\newblock Lecture Notes, Univ. Bayreuth, 2014.

\bibitem{degond7}
Andrew~V. Sutherland.
\newblock Genus 3 hyperelliptic curves.
\newblock \url{https://math.mit.edu/~drew/genus3curves.html}.

\end{thebibliography}
\end{document}